\documentclass{article}%
\usepackage{graphicx}
\usepackage{amsmath}
\usepackage{amsfonts}
\usepackage{amssymb}%
\providecommand{\U}[1]{\protect\rule{.1in}{.1in}}
\providecommand{\U}[1]{\protect\rule{.1in}{.1in}}
\providecommand{\U}[1]{\protect\rule{.1in}{.1in}}
\providecommand{\U}[1]{\protect\rule{.1in}{.1in}}
\providecommand{\U}[1]{\protect\rule{.1in}{.1in}}
\providecommand{\U}[1]{\protect\rule{.1in}{.1in}}
\providecommand{\U}[1]{\protect\rule{.1in}{.1in}}
\providecommand{\U}[1]{\protect\rule{.1in}{.1in}}
\providecommand{\U}[1]{\protect\rule{.1in}{.1in}}
\providecommand{\U}[1]{\protect\rule{.1in}{.1in}}
\providecommand{\U}[1]{\protect\rule{.1in}{.1in}}
\providecommand{\U}[1]{\protect\rule{.1in}{.1in}}
\providecommand{\U}[1]{\protect\rule{.1in}{.1in}}
\providecommand{\U}[1]{\protect\rule{.1in}{.1in}}
\providecommand{\U}[1]{\protect\rule{.1in}{.1in}}
\providecommand{\U}[1]{\protect\rule{.1in}{.1in}}
\providecommand{\U}[1]{\protect\rule{.1in}{.1in}}
\providecommand{\U}[1]{\protect\rule{.1in}{.1in}}
\providecommand{\U}[1]{\protect\rule{.1in}{.1in}}
\providecommand{\U}[1]{\protect\rule{.1in}{.1in}}
\providecommand{\U}[1]{\protect\rule{.1in}{.1in}}
\providecommand{\U}[1]{\protect\rule{.1in}{.1in}}
\providecommand{\U}[1]{\protect\rule{.1in}{.1in}}
\providecommand{\U}[1]{\protect\rule{.1in}{.1in}}
\providecommand{\U}[1]{\protect\rule{.1in}{.1in}}
\providecommand{\U}[1]{\protect\rule{.1in}{.1in}}
\providecommand{\U}[1]{\protect\rule{.1in}{.1in}}
\providecommand{\U}[1]{\protect\rule{.1in}{.1in}}
\providecommand{\U}[1]{\protect\rule{.1in}{.1in}}

\newtheorem{theorem}{Theorem}
{}

\newtheorem{definition}{Definition}
\newtheorem{example}{Example}

{}

\begin{document}

\title{A Brief Explanation of the Spectral Expansion Method for Non-Self-Adjoint
Differential Operators with Periodic Coefficients}
\author{O. A. Veliev\\{\small \ Dogus University, }\\{\small Esenkent 34755, \ Istanbul, Turkey.}\\\ {\small e-mail: oveliev@dogus.edu.tr}}
\date{}
\maketitle

\begin{abstract}
In this paper, we briefly explain the spectral expansion problem for
differential operators defined on the entire real line, generated by a
differential expression with periodic, complex-valued coefficients.



\end{abstract}

The aim is to provide a concise overview of the spectral expansion problem for
a non-self-adjoint differential operator $T(n)$ generated in $L_{2}%
(-\infty,\infty)$ by differential expression
\begin{equation}
y^{(n)}(x)+%
%TCIMACRO{\tsum \limits_{v=2}^{n}}%
%BeginExpansion
{\textstyle\sum\limits_{v=2}^{n}}
%EndExpansion
p_{v}(x)y^{(n-v)}(x) \tag{1}%
\end{equation}
with periodic complex-valued coefficients $p_{2},p_{3},...p_{n}$, where
$n\geq2.$ I will explain the method of constructing the spectral expansion and
the obtained results. Note that the main direct problems of differential
operators involve analyzing their spectrum and constructing the spectral
expansion. For self-adjoint operators in a Hilbert space, there is a
well-developed general theory, including the spectral theorem, which can be
applied to obtain the spectral expansion of the self-adjoint differential
operator. However, no general spectral expansion theorem exists for
non-self-adjoint operators. Moreover, well-known general theory of spectral
operators (see[3]) cannot be applied to the operator $T(n).$ Therefore, it was
necessary to develop a new method for constructing the spectral expansion in
this case.

For simplicity, we focus on the spectral theory of the one-dimensional
Schr\"{o}dinger operator $L(q)$ generated in $L_{2}(-\infty,\infty)$ by the
differential expression
\begin{equation}
l(y)=-y^{\prime\prime}+qy, \tag{2}%
\end{equation}
where $q$ is a complex-valued, $1$-periodic potential that is integrable on
$[0,1].$ However, the discussed method of constructing the spectral expansion
can also be applied to the operator $T(n).$ It is well-known that (see [7, 8])
the spectrum $\sigma(L(q))$ of the operator $L(q)$ is the union of the spectra
$\sigma(L_{t}(q))$ of the operators $L_{t}(q)$, for $t\in(-\pi,\pi]$,
generated in $L_{2}\left[  0,1\right]  $ by (2) and the boundary conditions
\begin{equation}
y\left(  1\right)  =e^{it}y\left(  0\right)  ,\text{ }y^{^{\prime}}\left(
1\right)  =e^{it}y^{^{\prime}}\left(  0\right)  . \tag{3}%
\end{equation}
The spectra $\sigma(L_{t})$ of the operators $L_{t}(q)$ consist of the
eigenvalues called the Bloch eigenvalues of $L(q)$. The eigenvalues
$\lambda_{k}(t)$ for $k\in\mathbb{Z}$ of $L_{t}(q)$ are known as Bloch
eigenvalues. They can be labeled so that for each $k\in\mathbb{Z}$ the
function $\lambda_{k}$ continuously depend on $t$ (see [18], Section 2.3).

First of all, let us note that the operator $L(q)$ is, in general, not a
spectral operator. Moreover, neither the smallness nor the smoothness of the
potential $q$ guarantees the spectrality of $L(q)$. Establishing spectrality
requires proving that the spectral projections of $L(q)$, corresponding to
subsets of the spectrum, are uniformly bounded. I proved that (see [10] and
Section 2.4 of [18]) if the geometric multiplicity of the multiple eigenvalue
$\lambda_{k}(t)$ is $1$ (it is clear that the geometric multiplicity of the
eigenvalues of $L_{t}(q)$ for $t\neq0,\pi$ is always $1$ (see, for example,
Proposition 2.2.4 of [18]), then the projections corresponding to the parts of
the spectrum lying in a neighborhood of this eigenvalue are not uniformly
bounded, that is, $\lambda_{k}(t)$ is a spectral singularity. Since this is a
widespread case for non-self-adjoint operators, $L(q)$ is, in general, not a
spectral operator.

The theory of spectral operators was developed to extend the spectral
expansion theorem from the class of self-adjoint operators to a wider class of
operators. However, the following examples demonstrate that, in some important
cases, the theory of spectral operators does not extend the spectral expansion
theorem from self-adjoint operators to a broader class of operators; that is,
the desired extension does not arise in these cases.

\begin{example}
Consider the case
\begin{equation}
q(x)=ae^{i2\pi x}+be^{-i2\pi x} \tag{4}%
\end{equation}
of Mathier-Schrodinger operator. I proved (see [15] and Chapter 4 of [18])
that if $ab\in\mathbb{R}$, then $L(q)$ is a spectral operator if and only it
is self adjoint. Thus, in this case, there are no spectral operators that are
not self-adjoint.
\end{example}

\begin{example}
Let $q$ be PT-symmetric periodic optical potential
\[
q(x)=4\cos^{2}x+4iV\sin2x=2+(1+2V)e^{i2x}+(1-2V)e^{-i2x},\text{ }V\geq0.
\]
I proved (see [11,13,16] and Chapter 5 of [18]) that the operator $L(q)$ is a
spectral operator if and only if $V=0$, that is, $q(x)=2+2\cos2x$ and $L(q)$
is a self-adjoint operator. Thus, in this important case as well, there are no
spectral operators that are not self-adjoint.
\end{example}

Let us redenote the operators $L(q)$ and $L_{t}(q)$ by $H(a,b)$ and
$H_{t}(a,b),$ respectively, when $q(x)$ has the form (4) and discuss the
spectrality problem for general Mathier-Schrodinger operator $H(a,b).$ Djakov
and Mityagin [1,2] proved that the root functions of the operators
$H_{0}(a,b)$ or $H_{\pi}(a,b)$ generated in $L_{2}[0,1]$ by periodic or
antiperiodic boundary conditions (see (3) for $t=0$ or $t=\pi$) form a Riesz
basis if and only if $\mid a\mid=\mid b\mid$ (see Theorem 7 in [2]). Moreover,
in [2] (see page 539) it was noted that Theorem 7 in [2] with Remark 8.10 in
[6] imply that $\ H(a,b)$ is not a spectral operator, if $|a|\neq\left\vert
b\right\vert .$ For another proof of this statement, see Section 4.3 in [18]
(Proposition 4.3.1). It is important to emphasize that the spectrality of
$H_{0}(a,b)$ and $H_{\pi}(a,b)$ (the periodic and antiperiodic cases), and
even the spectrality of all $H_{t}(a,b)$ for $t\in(-\pi,\pi],$ does not imply
the spectrality of $H(a,b)$ (see Remark 8.10 in [6]). Therefore, the condition
$\mid a\mid=\mid b\mid$ is not sufficient for the spectrality of $H(a,b)$.
Indeed, in [15] (see Theorem 11 in [15] or Theorem 4.1.1, page 190, in [18]),
I proved that $H(a,b)$ is an asymptotically spectral operator if and only if
$\mid a\mid=\mid b\mid$ and
\begin{equation}
\text{ }\inf_{q,p\in\mathbb{N}}\{\mid q\alpha-(2p-1)\mid\}\neq0,\tag{5}%
\end{equation}
where $\alpha=\pi^{-1}\arg(ab)$, $\mathbb{N=}\left\{  1,2,...,\right\}  .$ For
example, if $\alpha=\frac{m}{q},$ where $m$ is an odd integer and $m,q$ are
coprime, then $H(a,b)$ is not even an asymptotically spectral operator. In
other word, even when $\mid a\mid=\mid b\mid$ if condition (5) fails, then
$H(a,b)$ is not a spectral operator. Thus Theorem 7 of [2] provide a necessary
condition for the spectrality of $H(a,b)$, but not a sufficient one, that is,
does not characterize the spectrality of $H(a,b)$ fully. That is why, in the
first version of this paper (v1), I discussed only the special cases (Example
1 and Example 2), for which we have necessary and sufficient condition for
spectrality. Moreover, in Remark 4.3.5 of [18], I noted that a detailed
investigation of $H_{0}(a,b)$ and $H_{\pi}(a,b)$ cannot provide a sufficient
condition for the spectrality of $\ H(a,b)$. This is because multiple
eigenvalues of $H_{t}(a,b)$ for $t\neq0,\pi$ become spectral singularities of
$H(a,b)$, which in turn implies that $H(a,b)$ is not a spectral operator. Such
multiple eigenvalues may indeed arise for $t\neq0,\pi$ (see Theorem 4.3.1 in [18]).

Thus, it is natural to conclude that finding explicit conditions on $q$ that
guarantee the spectrality of $L(q)$ is a complex and generally ineffective
problem and the spectrality of $L(q)$ is a very rare occurrence. However, in
[6], Gesztesy and Tkachenko established two versions of a criterion for $L(q)$
to be a spectral operator in the sense of Dunford---one analytic and one
geometric. The analytic version was formulated in terms of the solutions of
the Hill equation, while the geometric version was expressed through algebraic
and geometric properties of the spectra of the periodic/antiperiodic and
Dirichlet boundary value problem. The problem of explicitly characterizing the
potentials $q$ for which the Schr\"{o}dinger operators $L(q)$ are spectral
operators has remained open for approximately 65 years. In [9], I found
explicit conditions on the potential $q$ under which $L(q)$ is an
asymptotically spectral operator, using asymptotic formulas. However, since
these formulas do not provide information about the small spectral
singularities, my method also does not yield any explicit conditions on $q$
for the spectrality of $L(q).$

Now, let us return to discussing the spectral expansion problem for the
Schr\"{o}dinger operators $L(q).$ Gelfand's method for the spectral expansion
of self-adjoint differential operators with periodic coefficients is based on
the following (see [5]). For every $f\in L_{2}(-\infty,\infty)$ there exists%
\[
f_{t}(x)=\sum_{k=-\infty}^{\infty}f(x+k)e^{-ikt}%
\]
such that
\begin{equation}
f(x)=\frac{1}{2\pi}\int\limits_{-\pi}^{\pi}f_{t}(x)dt,\text{ }\int
\limits_{-\infty}^{\infty}\left\vert f(x)\right\vert ^{2}dx=\frac{1}{2\pi}%
\int\limits_{-\pi}^{\pi}\int\limits_{0}^{1}\left\vert f_{t}(x)\right\vert
^{2}dxdt. \tag{6}%
\end{equation}
In the self-adjoint case, Parseval's equality for the operator $L_{t}(q)$
implies that%
\begin{equation}
\int\limits_{0}^{1}\left\vert f_{t}(x)\right\vert ^{2}dx=\sum_{n\in\mathbb{Z}%
}\left\vert a_{n}(t)\right\vert ^{2},\text{ }\int\limits_{-\infty}^{\infty
}\left\vert f(x)\right\vert ^{2}dx=\frac{1}{2\pi}\int\limits_{-\pi}^{\pi}%
\sum_{n\in\mathbb{Z}}\left\vert a_{n}(t)\right\vert ^{2}dt, \tag{7}%
\end{equation}
where
\[
a_{n}(t)=\int\nolimits_{0}^{1}f_{t}(x)\overline{\Psi_{n,t}(x)}dx
\]
and $\Psi_{n,t}$ is the normalized eigenfunctions of $L_{t}(q).$

It is easy to see that series in (7) can be integrated term by term. This
yields the decomposition
\begin{equation}
\int\limits_{-\infty}^{\infty}\left\vert f(x)\right\vert ^{2}dx=\frac{1}{2\pi
}\sum_{n\in\mathbb{Z}}\int\nolimits_{-\pi}^{\pi}\left\vert a_{n}(t)\right\vert
^{2}dt,\text{ }f=\frac{1}{2\pi}\sum_{n\in\mathbb{Z}}\int\nolimits_{-\pi}^{\pi
}a_{n}(t)\Psi_{n,t}dt. \tag{8}%
\end{equation}
Thus, in the self-adjoint case, the existence of Parseval's equality readily
implies the spectral expansion (8). However, in the non-self-adjoint case, no
such Parseval equality holds. The absence of Parseval's equality for the
non-self-adjoint operator $L_{t}(q)$ prevents the application of Gelfand's
elegant method to construct the spectral expansion for the non-self-adjoint
operator $L(q).$

Since the $t$-periodic boundary conditions are strongly regular for all
$t\neq0,\pi$, we can employ the Riesz basis property of the root functions of
$L_{t}(q)$ for $t\in(-\pi,\pi]\backslash\left\{  0,\pi\right\}  .$ Moreover,
for almost all $t\in(-\pi,\pi]$, all eigenvalues of the operators $L_{t}(q)$
are simple. Therefore, the system $\{\Psi_{n,t}:n\in\mathbb{Z}\}$ of
normalized eigenfunctions of $L_{t}(q)$ is a Reisz basis of $L_{2}[0,1]$ and
thus we have the decomposition
\begin{equation}
f_{t}=\sum_{n\in\mathbb{Z}}a_{n}(t)\Psi_{n,t} \tag{9}%
\end{equation}
of $f_{t}$ by the basis $\{\Psi_{n,t}:n\in\mathbb{Z}\},$ for almost all
$t\in(-\pi,\pi],$ where $a_{n}(t)=(f_{t},X_{n,t})$ and $\left\{  X_{n,t}%
:n\in\mathbb{Z}\right\}  $ is the biorthogonal system, that is,
\[
X_{n,t}=\frac{1}{\alpha_{n}(t)}\Psi_{n,t}^{\ast},\text{ }\alpha_{n}(t)=\left(
\Psi_{n,t},\Psi_{n,t}^{\ast}\right)  ,
\]
and $\Psi_{n,t}^{\ast}$ is the normalized eigenfunctions of $(L_{t}(q))^{\ast
}$ .

Using (9) in (6), we get
\begin{equation}
f=\frac{1}{2\pi}\int\limits_{(-\pi,\pi]}\sum_{n\in\mathbb{Z}}a_{n}%
(t)\Psi_{n,t}dt, \tag{10}%
\end{equation}
where
\begin{equation}
a_{n}(t)\Psi_{n,t}=\frac{1}{\alpha_{n}(t)}(f_{t},\Psi_{n,t}^{\ast})\Psi_{n,t}.
\tag{11}%
\end{equation}
To obtain the spectral expansion in terms of $t$ from (10), it is necessary to
justify the term-by-term integration. This has been a challenging and
intricate problem since the 1950s. Moreover, the term-by term integration in
(10) is, in general, not possible, because the expression (11) is not
integrable over $(-\pi,\pi]$ for certain values of $n.$

Thus, the first step is to examine the integrability of (11). Since the
integrability of $a_{n}(t)\Psi_{n,t}$ depends on the integrability of
$\frac{1}{\alpha_{n}(t)}$, we introduce the following notions, defined
independently of the choice of $f$ , which will be used in the construction of
the spectral expansion.

\begin{definition}
We say that a point $\lambda_{0}\in\sigma(L_{t_{0}}(q))\subset\sigma(L(q))$ is
an essential spectral singularity (ESS) of the operator $L$ if there exists
$n\in\mathbb{Z}$ such that $\lambda_{0}=:\lambda_{n}(t_{0})$ and for each
$\varepsilon$ the function $\frac{1}{\alpha_{n}}$ is not integrable on
$(t_{0}-\varepsilon,t_{0}+\varepsilon)$. Then the quasimomentum $t_{0}$ is
said to be singular quasimomentum (SQ).
\end{definition}

\begin{definition}
We say that the operator $L$ has ESS at infinity if there exist a sequence of
integers $k_{s}$ and a sequence of closed subsets $I(s)$ of $(-\pi,\pi]$ such
that
\[
\lim_{s\rightarrow\infty}\int\nolimits_{I(s)}\left\vert \alpha_{k_{s}%
}(t)\right\vert ^{-1}dt=\infty.
\]

\end{definition}

Using these notions, we proved the following (see [12] and Chapter 2 of [18]):

\begin{theorem}
We have the elegant spectral decompositions
\begin{equation}
f(x)=\frac{1}{2\pi}\sum_{k\in\mathbb{Z}}\int_{0}^{2\pi}a_{k}(t)\Psi
_{k,t}(x)dt=\int\limits_{\sigma(L(q))}\Phi(x,\lambda)d\lambda\tag{12}%
\end{equation}
if and only if $L(q)$ has no ESS and ESS at infinity, where
\[
\Phi(x,\lambda)=(\Phi_{t}(x,\lambda)F_{-}(\lambda,f)+\Phi_{-t}(x,\lambda
)F_{+}(\lambda,f))\frac{1}{\varphi(1,\lambda)p(\lambda)},
\]%
\[
\Phi_{t}(x,\lambda)=\varphi(1,\lambda)\theta(x,\lambda)+(e^{it}-\theta
(1,\lambda))\varphi(x,\lambda),
\]
for $\lambda\in\sigma(L_{t}(q)),$%
\[
\text{ }p(\lambda)=\sqrt{4-F^{2}(\lambda)},\text{ }F_{\pm}(\lambda
,f)=\int_{\mathbb{R}}f(x)\Phi_{\pm t}(x,\lambda)dx,
\]
$F(\lambda)=\varphi^{\prime}(1,\lambda)+\theta(1,\lambda),$ $\varphi
(x,\lambda)$ and $\theta(x,\lambda)$ are the solutions of
\[
-y^{\prime\prime}(x)+q(x)y(x)=\lambda y
\]
satisfying the conditions $\varphi(0,\lambda)=0,$ $\varphi^{\prime}%
(0,\lambda)=1$ and $\theta(0,\lambda)=1,$ $\theta^{\prime}(0,\lambda)=0.$
\end{theorem}

In [15] and Chapter 4 of [18], we proved that if $q(x)$ has the form (4) and
$0<\left\vert ab\right\vert <\frac{16}{9}\pi^{4}$, then $L(q)$ has neither ESS
nor an ESS at infinity, and (12) holds. Moreover, in the case of the optical
potential (see Example 2), we proved that there exists a sequence of critical
points $V_{1}=1/2,V_{2}=0.888437...,...,$ $($with $V_{k}\rightarrow\infty$ as
$k\rightarrow\infty)$ such that if $V\neq V_{k},$ then $L(q)$ has infinitely
many spectral singularity, has neither ESS nor an ESS at infinity, and (12)
remains valid (see [16] and Chapter 5 of [18]).

Now let us discuss the case when $L(q)$ has ESS or ESS at infinity. In [12]
and Chapter 2 of [18], we have established that only the periodic eigenvalues
and the antiperiodic eigenvalues may become ESS. Let $\Lambda_{1}%
(0),\Lambda_{2}(0),...$ and $\Lambda_{1}(\pi),\Lambda_{2}(\pi),...$ denote,
respectively, the periodic and antiperiodic eigenvalues that are ESS. For each
$\Lambda_{j}(0)$ and $\Lambda_{j}(\pi)$ introduce the notation
\[
\mathbb{S}(\Lambda_{j}(0)):=\left\{  n\in\mathbb{Z}:\lambda_{n}(0)=\Lambda
_{j}(0)\right\}  ,\text{ }\mathbb{S}(j,0):=\left\{  a_{k}(t)\Psi_{k,t}%
:k\in\mathbb{S}(\Lambda_{j}(0))\right\}
\]
and
\[
\mathbb{S}(\Lambda_{j}(\pi)):=\left\{  n\in\mathbb{Z}:\lambda_{n}(\pi
)=\Lambda_{j}(\pi)\right\}  ,\text{ }\mathbb{S}(j,\pi):=\left\{  a_{k}%
(t)\Psi_{k,t}:k\in\mathbb{S}(\Lambda_{j}(\pi))\right\}
\]
We have proved that every element of the sets $\mathbb{S}(j,0)$ and
$\mathbb{S}(j,\pi)$ is not integrable on $[-h,h]$ and $[\pi-h,\pi+h]$
respectively, where $0<h<1/15\pi.$ However, the sum of the elements of the
sets $\mathbb{S}(j,0)$ and $\mathbb{S}(j,\pi)$ is integrable on $[-h,h]$ and
$[\pi-h,\pi+h]$ respectively. Moreover, the following equalities are
satisfied
\begin{equation}
\int\limits_{\lbrack-h,h]}\sum\limits_{k\in\mathbb{S}(\Lambda_{j}(0))}%
a_{k}(t)\Psi_{k,t}dt=\lim_{\delta\rightarrow0}\sum\limits_{k\in\mathbb{S}%
(\Lambda_{j}(0))}\int\limits_{[-h,h]\backslash\lbrack-\delta,\delta]}%
a_{k}(t)\Psi_{k,t}dt \tag{13}%
\end{equation}
and%
\begin{equation}
\int\limits_{\lbrack\pi-h,\pi+h]}\sum\limits_{k\in\mathbb{S}(\Lambda_{j}%
(\pi))}a_{k}(t)\Psi_{k,t}dt=\lim_{\delta\rightarrow0}\sum\limits_{k\in
\mathbb{S}(\Lambda_{j}(\pi))}\int\limits_{[\pi-h,\pi+h]\backslash\lbrack
\pi-\delta,\pi+\delta]}a_{k}(t)\Psi_{k,t}dt. \tag{14}%
\end{equation}
Therefore, the elements of the set $\mathbb{S}(j,0)$ must be grouped together.
Similarly, the elements of the set $\mathbb{S}(j,\pi)$ must also be grouped together.

That is why we rewrite (6) in the following form:
\begin{equation}
f=\frac{1}{2\pi}\left(  \int\limits_{B(h)}f_{t}dt+\int\limits_{(-h,h)}%
f_{t}dt+\int\limits_{(\pi-h,\pi+h)}f_{t}dt\right)  , \tag{15}%
\end{equation}
where $B(h)=[-\pi+h,-h]\cup\lbrack h,\pi-h]$ and prove that
\begin{equation}
\int\limits_{B(h)}f_{t}dt=\sum\limits_{k\in\mathbb{Z}}\int\limits_{B(h)}%
a_{k}(t)\Psi_{k,t}dt, \tag{16}%
\end{equation}%
\begin{equation}
\int\limits_{(-h,h)}f_{t}dt=\sum_{j}\int\limits_{(-h,h)}\sum\limits_{k\in
\mathbb{S}(\Lambda_{j}(0))}a_{k}(t)\Psi_{k,t}dt+\sum_{k\in\mathbb{Z}%
\backslash\mathbb{S}(0)}\int\limits_{(-h,h)}a_{k}(t)\Psi_{k,t}dt, \tag{17}%
\end{equation}
and
\begin{equation}
\int\limits_{(\pi-h,\pi+h)}f_{t}dt=\sum_{j}\int\limits_{(\pi-h,\pi+h)}%
\sum\limits_{k\in\mathbb{S}(\Lambda_{j}(\pi))}a_{k}(t)\Psi_{k,t}dt+\sum
_{k\in\mathbb{Z}\backslash\mathbb{S}(\pi)}\int\limits_{(-h,h)}a_{k}%
(t)\Psi_{k,t}dt, \tag{18}%
\end{equation}
where $\mathbb{S}(0)=\cup_{j}(\mathbb{S}(\Lambda_{j}(0)))$ and $\mathbb{S}%
(\pi)=\cup_{j}(\mathbb{S}(\Lambda_{j}(\pi)))$ and equalities (13) and (14) for
the first integrals in the right-side of (17) and (18) are satisfied (see [12]
and Chapter 2 of [18]). Equalities (15)-(18) is the spectral expansion in term
of the quasimomentum $t\in(-\pi,\pi].$

To write the spectral expansion in term of the spectral parameter $\lambda
\in\sigma(L(q)),$ introduce the notation:
\[
\gamma(j,0,h)=\bigcup\limits_{n\in\mathbb{S}(\Lambda_{j}(0))}\left\{
\lambda_{n}(t):t\in(-h,h)\right\}  ,
\]%
\[
\gamma(j,\pi,h)=\bigcup\limits_{n\in\mathbb{S}(\Lambda_{j}(\pi))}\left\{
\lambda_{n}(t):t\in(\pi-h,\pi+h)\right\}  .
\]
and
\[
\gamma(h)=\left(  \cup_{j}\gamma(j,0,h)\right)  \cup(\cup_{j}\gamma
(j,\pi,h))).
\]
From (15)-(18), by changing the variables, we obtain the following spectral
expansion (see [12] and Chapter 2 of [18]):

\begin{theorem}
\textbf{(Spectral Expansion Theorem)}. The following spectral expansion
holds:
\begin{align*}
f(x)  &  =\int\limits_{\sigma(L)\backslash\gamma(h),}\Phi(x,\lambda
)d\lambda+\sum\limits_{j}p.v.\int\limits_{\gamma(j,0,h)}\Phi(x,\lambda
)d\lambda+\\
&  \sum\limits_{j}p.v.\int\limits_{\gamma(j,\pi,h)}\Phi(x,\lambda)d\lambda,,
\end{align*}
where the principal value integrals over $\gamma(j,0,h)$ and $\gamma(j,\pi,h)$
are defined as follows:%
\[
p.v.\int\limits_{\gamma(j,0,h)}\Phi(x,\lambda)d\lambda=\lim_{\delta
\rightarrow0}\int\limits_{\gamma(j,0,h)\backslash\gamma(j,0,\delta),}%
\Phi(x,\lambda)d\lambda
\]
and
\[
p.v.\int\limits_{\gamma(j,\pi,h)}\Phi(x,\lambda)d\lambda=\lim_{\delta
\rightarrow0}\int\limits_{\gamma(j,\pi,h)\backslash\gamma(j,\pi,\delta),}%
\Phi(x,\lambda)d\lambda.
\]

\end{theorem}

We now illustrate the spectral expansion of $L(q)$ in the case of an optical
potential. As was noted above, if $V\neq V_{k},$ then (12) remains valid. The
case, $V=V_{1}=1/2,$ that is, $q(x)=2+2e^{i2x}$ is the special potential,
which admits explicit special solutions and a special spectral expansion that
cannot be applied to other potentials (see [4]). If $V=V_{2}$ then there
exists one ESS $\lambda=\lambda_{1}(0)=\lambda_{2}(0)$, infinitely many
spectral singularity and we have the following spectral expansion in term of
$t:$
\[
f=%
%TCIMACRO{\tint \nolimits_{(-\pi,\pi]}}%
%BeginExpansion
{\textstyle\int\nolimits_{(-\pi,\pi]}}
%EndExpansion
(a_{1}(t)\Psi_{1,t}+a_{2}(t)\Psi_{2,t})dt+%
%TCIMACRO{\tsum _{k>2}}%
%BeginExpansion
{\textstyle\sum_{k>2}}
%EndExpansion%
%TCIMACRO{\tint \nolimits_{(-\pi,\pi]}}%
%BeginExpansion
{\textstyle\int\nolimits_{(-\pi,\pi]}}
%EndExpansion
a_{k}(t)\Psi_{k,t}dt.
\]
(see [16] and Chapter 5 of [18]). Similarly, if $V=V_{3}$ then there exists
one ESS $\lambda=\lambda_{3}(0)=\lambda_{4}(0)$, infinitely many spectral
singularity and the spectral expansion has the form
\[
f=%
%TCIMACRO{\tint \nolimits_{(-\pi,\pi]}}%
%BeginExpansion
{\textstyle\int\nolimits_{(-\pi,\pi]}}
%EndExpansion
(a_{3}(t)\Psi_{3,t}+a_{4}(t)\Psi_{4,t})dt+\sum_{k\neq3,4}%
%TCIMACRO{\tint \nolimits_{(-\pi,\pi]}}%
%BeginExpansion
{\textstyle\int\nolimits_{(-\pi,\pi]}}
%EndExpansion
a_{k}(t)\Psi_{k,t}dt.
\]

By the same method, I constructed the spectral expansion for high order
differential operator $T(n)$ generated in $L_{2}(-\infty,\infty)$ by (1) (see
[14] and Chapter 7 of [18]). Moreover, by developing this method, we
constructed the spectral expansion for the differential operator $T(n,m)$
generated in the space $L_{2}^{m}(-\infty,\infty)$ of vector-valued functions
by the differential expression
\[
l(y)=y^{(n)}(x)+P_{2}\left(  x\right)  y^{(n-2)}(x)+P_{3}\left(  x\right)
y^{(n-3)}(x)+...+P_{n}(x)y,
\]
if the eigenvalues $\mu_{1},\mu_{2},...,\mu_{m}$ of the matrix $\int_{0}%
^{1}P_{2}\left(  x\right)  dx$ are simple, where $n\geq2$ and $P_{\mathbb{\nu
}}$ for $\mathbb{\nu}=2,3,...n$ are the $m\times m$ matrix with the
complex-valued periodic entries $\left(  p_{\mathbb{\nu},i,j}\right)  $ (see
[17] and Chapters 6 and 7 of [18]).

\end{document}